\documentclass[12pt]{article}
\usepackage{amsmath, amsthm, amsfonts}
\usepackage[T1]{fontenc}
\usepackage{lmodern}
\usepackage[utf8]{inputenc}
\usepackage{ragged2e}
\usepackage{mathtools}
\usepackage{lipsum}
\usepackage[top=1in, bottom=1in, left=1in, right=1in]{geometry}
\usepackage{microtype}
\usepackage{amssymb}
\usepackage[hidelinks]{hyperref}
\usepackage{listings}
\numberwithin{equation}{section}
\usepackage{booktabs}
\usepackage{hyperref}
\usepackage[nameinlink,noabbrev]{cleveref}

\newtheorem{theorem}{Theorem}[section]
\newtheorem{lemma}[theorem]{Lemma}
\newtheorem{proposition}[theorem]{Proposition}
\newtheorem{corollary}[theorem]{Corollary}

\begin{document}

\vspace{100pt}

\begin{center}
    \textbf{\LARGE Bounded Ratios of Best Approximations: Sharp Bounds and Spectra of Diophantine Exponents}
\end{center}

\begin{center}
    \textbf{\Large Maksim Lavrov}\\[4pt]
    Faculty of Mechanics and Mathematics,\\
    Lomonosov Moscow State University
\end{center}

\begin{abstract}
In a recent paper, Marnat, Moshchevitin, and Schleischitz studied
Diophantine exponents of totally irrational matrices whose sequences
of best approximations satisfy boundedness conditions on the
successive height ratios $X_{i+1}/X_i$ and error ratios
$L_i/L_{i+1}$. They posed several problems concerning the possible
values of the ordinary and uniform Diophantine exponents under these
conditions.

Let $\xi$ be a totally irrational $n\times m$ matrix. We solve Problems~1--3 posed in that paper. In particular, assuming only that the error ratios $L_i/L_{i+1}$ are bounded, we prove the sharp upper bounds
\[
    \omega(\xi) \leq m
    \qquad\text{and}\qquad
    \hat{\omega}(\xi) \leq \frac{m^2}{m+1}
\]
for $m,n\geq 2$. We then determine the corresponding spectra, both
with and without the additional assumption that the height ratios
$X_{i+1}/X_i$ are bounded. Our proofs are based on a uniform estimate
for blocks of consecutive best approximations contained in a fixed
rational subspace, together with a geometric analysis of transitions
between such subspaces.
\end{abstract}

\section{Introduction}

Let $m,n\geq 1$, $d=m+n$, and $\xi\in \mathbb{R}^{n\times m}$. For $z=(x,y)\in \mathbb{R}^m \oplus \mathbb{R}^n$ define
\begin{equation}
    X(z):=\|x\|_2, \quad L(z):=\|\xi x-y\|_2, \quad \mathcal{L}_\xi :=\{(x,y): \xi x-y=0\}.
\end{equation}
We consider two conditions on the matrix $\xi$:
\begin{equation}\label{cond1}
    \mathcal{L}_\xi \cap \mathbb{Z}^d=\{0\},
\end{equation}
and
\begin{equation}\label{cond2}
    \mathcal{L}_\xi \ \text{is not contained in any proper rational subspace of} \ \mathbb{R}^d.
\end{equation}
A matrix is called totally irrational if it satisfies \eqref{cond1} and \eqref{cond2}. A subspace is called rational if it is spanned by rational vectors.

Choose a sequence of best approximations $z_i=(x_i,y_i)\in \mathbb{Z}^d$, and write $X_i:=X(z_i)$, $L_i:=L(z_i)$. We assume that \eqref{cond1} holds. Then
\[
X_1<X_2<\dots \to \infty, \quad L_1>L_2>\dots \to 0.
\]
The defining property of this sequence is
\begin{equation}\label{defining property}
    0\neq z\in\mathbb{Z}^d, \quad X(z)<X_{i+1} \implies L(z)\geq L_i,
\end{equation}
for all sufficiently large $i$.

Following \cite{MMS}, we introduce two conditions. There exists a constant $M_1>0$ such that for all sufficiently large $i$
\begin{equation*}\label{cond B}\tag{B}
    \frac{X_{i+1}}{X_i}\leq M_1.
\end{equation*}
There exists a constant $M_2>0$ such that for all sufficiently large $i$
\begin{equation*}\label{cond C}\tag{C}
    \frac{L_i}{L_{i+1}}\leq M_2.
\end{equation*}

Throughout the paper, best approximations are taken with respect to
the Euclidean norms. For totally irrational matrices, the boundedness
conditions \eqref{cond B} and \eqref{cond C} are independent of the
choice of norms by \cite[Theorem~2.1]{MMS}. When only \eqref{cond1}
is assumed, we work directly with the Euclidean best approximations.

In \cite{MMS}, the authors pose four problems. The main goal of this paper is to give solutions to three of them.

\noindent\textbf{Problem 1:}
For $m \geq 3$ and $n \geq 1$, find an upper bound for the
ordinary exponent $\omega$ of totally irrational matrices
$\xi \in \mathbb{R}^{n \times m}$ satisfying \eqref{cond C}.

The authors conjecture that the optimal bound is $m$. They prove that this is the case for $n=1$, and $n\geq 2, m=1$, and $n\geq 2, m=2$. For $m\geq3$ and $n\geq2$, they prove that $\omega(\xi)<\infty$ for every totally irrational matrix satisfying \eqref{cond C} \cite[Theorems~2.2 and~2.8]{MMS}.

\noindent\textbf{Problem 2:}
For $m\geq 2$ and $n\geq 2$, find an upper bound for the uniform exponent $\hat{\omega}$ of totally irrational matrices $\xi\in\mathbb{R}^{n\times m}$ satisfying \eqref{cond C}.

The authors conjecture that the optimal bound is $\frac{m^2}{m+1}$. Note that in the original paper the problem is stated for $n\geq 1$. However, in the case $n=1$ condition \eqref{cond C} implies that $\xi$ is badly approximable and the problem is trivial, since $\omega(\xi)=\hat{\omega}(\xi)=m$. The optimality of this bound in the case $m=1$ is proven by Summerer in \cite{Summerer}. In \cite{MMS} the authors show that under condition \eqref{cond C} any uniform exponent from the interval $\left[\frac{m}{n},\frac{m^2}{m+1}\right]$ is admissible and prove that in the case where \eqref{cond B} also holds and $m=2$, $\hat{\omega}(\xi)\leq 4/3$ \cite[Theorems~2.9,2.10]{MMS}.

\noindent\textbf{Problem 3:} Are the upper bounds for ordinary and/or uniform exponents the same for totally irrational matrices $\xi\in\mathbb{R}^{n\times m}$ satisfying only \eqref{cond C} and totally irrational matrices $\xi\in\mathbb{R}^{n\times m}$ satisfying both \eqref{cond B} and \eqref{cond C}?

In the present paper we give solutions to all these problems. In particular, we prove the following theorems.

\begin{theorem}[Block lemma]\label{th:block lemma}
    Assume that \eqref{cond1} and \eqref{cond C} hold. Then there exist an index $a_0$ and a constant $C_*=C_*(m,n,M_2)>0$ such that for all $a\geq a_0$ the following holds. Suppose $V$ is a rational subspace of $\mathbb{R}^d$ of dimension $r\geq 2$, and
    \[
    z_a,\dots, z_b\in V, \quad a<b.
    \]
    Set $\mu=\min(m,r-1)$. Then
    \begin{equation}\label{block lemma ineq}
        L_iX_{i+1}^\mu \leq C_* L_jX_{j+1}^\mu \quad (a\leq i\leq j<b).
    \end{equation}
\end{theorem}

The block lemma yields the following solution to Problem~1.

\begin{theorem}\label{problem 1 sol}
    Under \eqref{cond1} and \eqref{cond C} there exists a constant $c_\xi>0$ such that
    \begin{equation}\label{theorem 2 ineq}
        L_iX_i^m\geq c_\xi \quad \text{for all sufficiently large}\ i.
    \end{equation}
    Consequently,
    \begin{equation}
        \omega(\xi)\leq m.
    \end{equation}
\end{theorem}

\Cref{problem 1 sol} requires only \eqref{cond1} and \eqref{cond C}, without assuming \eqref{cond2}. In particular, it applies to totally irrational matrices satisfying \eqref{cond C}. It is shown in \cite[Theorem~2.6(2)]{MMS} that every exponent from $[m/n,m]$ is admissible for totally irrational matrices when $n\geq 2$ and condition \eqref{cond C} holds (if $n=1$ the problem is trivial since $\xi$ becomes badly approximable). This implies that our bound is indeed optimal and solves the first problem.

As a further application of the block lemma, we obtain the following result.

\begin{theorem}\label{tail theorem}
    Let $m,n\geq 2$. Assume \eqref{cond1}, \eqref{cond2}, \eqref{cond C}. If $z_i$ lies in a fixed rational subspace $Q$ with $\dim Q\leq m+1$ for all sufficiently large $i$, then
    \begin{equation}\label{tail theorem ineq}
        L_iX_i^{m-1} \gg_{\xi,Q,M_2} 1, \quad \hat{\omega}(\xi)\leq \omega(\xi)\leq m-1. 
    \end{equation}
\end{theorem}

The next theorem solves the second problem.

\begin{theorem}\label{problem 2 sol}
    Let $m,n\geq 2$. Assume \eqref{cond1}, \eqref{cond2}, \eqref{cond C}. Then
    \[
    \hat{\omega}(\xi)\leq \frac{m^2}{m+1}.
    \]
\end{theorem}

Finally, we turn to Problem~3. For $m,n\geq 1$ we denote by $\mathcal{C}_{m,n}$ the set of all totally irrational matrices satisfying \eqref{cond C} and by $\mathcal{B}\mathcal{C}_{m,n}$ the set of all totally irrational matrices satisfying both \eqref{cond B} and \eqref{cond C}. For a class $\mathcal{A}$ and an exponent $\tau \in \{\omega,\hat{\omega}\}$ define
\[
\Sigma_\tau (\mathcal{A}):=\{ \tau(\xi):\xi\in\mathcal{A}\}.
\]
The following theorem answers Problem~3.

\begin{theorem}\label{problem 3 sol}
    The spectra are given in the following table.
    \begin{table}[ht]
        \centering
        \renewcommand{\arraystretch}{1.4}
        \begin{tabular}{@{}lcccc@{}}
        \toprule
        \textit{Dimensions}
        & $\Sigma_{\omega}(\mathcal{C}_{m,n})$
        & $\Sigma_{\omega}(\mathcal{BC}_{m,n})$
        & $\Sigma_{\hat{\omega}}(\mathcal{C}_{m,n})$
        & $\Sigma_{\hat{\omega}}(\mathcal{BC}_{m,n})$ \\
        \midrule
        $m,n \geq 2$
        & $\left[\dfrac{m}{n},m\right]$
        & $\left[\dfrac{m}{n},m\right]$
        & $\left[\dfrac{m}{n},\dfrac{m^{2}}{m+1}\right]$
        & $\left[\dfrac{m}{n},\dfrac{m^{2}}{m+1}\right]$ \\
        
        $n=1,\ m\geq 1$
        & $\{m\}$
        & $\{m\}$
        & $\{m\}$
        & $\{m\}$ \\
        
        $m=1,\ n\geq 2$
        & $\left[\dfrac{1}{n},1\right]$
        & $\left\{\dfrac{1}{n}\right\}$
        & $\left[\dfrac{1}{n},\dfrac{1}{2}\right]$
        & $\left\{\dfrac{1}{n}\right\}$ \\
        \bottomrule
        \end{tabular}
    \end{table}
\end{theorem}

We shall also use the following standard formulas:
\begin{equation}\label{eq:exp formulas}
    \omega(\xi)=\limsup_{i\to\infty} \frac{-\log L_i}{\log X_i},\quad \hat{\omega}(\xi)=\liminf_{i\to\infty}\frac{-\log L_i}{\log X_{i+1}}.
\end{equation}
Dirichlet's theorem and the definitions give $m/n\leq\hat{\omega}(\xi)\leq\omega(\xi)$.

\section{General theory}

In this section, we develop the tools used in the proofs of the main results. Throughout this section, we assume \eqref{cond1} and \eqref{cond C}.

\subsection{Properties of the first minimum}

Define
\[
E:=\mathbb{R}^m\oplus0,\quad F:=0\oplus\mathbb{R}^n,\quad \Phi(x,y):=(x,\xi x-y), \ x\in \mathbb{R}^m, \ y\in \mathbb{R}^n,
\]
and
\[
D_t(u,v):=(e^{-nt}u,e^{mt}v), \quad A_t:=D_t\Phi.
\]
Clearly all of these operators are invertible. On $E\oplus F$ we use the norm
\[
\|(u,v)\|_\square:=\max (\|u\|_2,\|v\|_2).
\]
For this norm
\[
\|(u,v)\|_2^2=\|u\|_2^2+\|v\|_2^2\leq 2\max(\|u\|_2^2,\|v\|_2^2)=2\|(u,v)\|_\square ^2.
\]
This implies
\begin{equation}\label{block and eucl norm comp}
    \|z\|_\square\leq \|z\|_2\leq \sqrt{2}\|z\|_\square.
\end{equation}

The logarithm of the first minimum of the full lattice $A_t\mathbb{Z}^d$ is
\begin{equation}\label{h1}
    h_1(t):=\log \min_{0\neq z\in \mathbb{Z}^d}\|A_tz\|_\square.
\end{equation}
Define
\begin{equation}\label{pi qi}
    p_i:=\frac{\log X_i-\log L_i}{d}, \quad q_i:=\frac{\log X_{i+1}-\log L_i}{d}.
\end{equation}
Then $p_i<q_i<p_{i+1}$, and the following lemma holds.

\begin{lemma}
    For sufficiently large $i$
    \begin{equation}\label{h1 formula}
    h_1(t)=
    \begin{cases}
    \log L_i + mt,      & p_i \leq t \leq q_i,\\
    \log X_{i+1} - nt,  & q_i \leq t \leq p_{i+1}.
    \end{cases}
    \end{equation}
\end{lemma}
\begin{proof}
    Fix $t\in[p_i,p_{i+1}]$. Let $\lambda_1(t)=e^{h_1(t)}$. By definition of $A_t$ and $\|\cdot\|_\square$ we have
    \[
    \|A_tz\|_\square=\max (e^{-nt}X(z),e^{mt}L(z)).
    \]
    Then
    \[
    \|A_tz_i\|_\square=e^{mt}L_i, \quad \|A_tz_{i+1}\|_\square=e^{-nt}X_{i+1}.
    \]
    The first equality holds since the inequality $e^{-nt}X_i\leq e^{mt}L_i$ is equivalent to $t\geq p_i$. And the second equality holds since $e^{-nt}X_{i+1}\geq e^{mt}L_{i+1}$ is equivalent to $t\leq p_{i+1}$. Therefore
    \[
    \lambda_1(t)\leq \min(e^{mt}L_i,e^{-nt}X_{i+1}).
    \]

    To prove the reverse inequality consider any vector $0\neq z\in\mathbb{Z}^d$. If $X(z)<X_{i+1}$, the best-approximation property \eqref{defining property} implies
    \[
    \|A_tz\|_\square\geq e^{mt}L(z)\geq e^{mt}L_i\geq \min(e^{mt}L_i,e^{-nt}X_{i+1}).
    \]
    And if $X(z)\geq X_{i+1}$, then
    \[
    \|A_tz\|_\square\geq e^{-nt}X(z)\geq e^{-nt}X_{i+1}\geq\min(e^{mt}L_i,e^{-nt}X_{i+1}).
    \]
    Hence
    \[
    \|A_tz\|_\square\geq \min(e^{mt}L_i,e^{-nt}X_{i+1}).
    \]
    Taking the minimum over all nonzero integer $z$ proves the reverse inequality. Thus
    \[
    \lambda_1(t)=\min(e^{mt}L_i,e^{-nt}X_{i+1}).
    \]

    Finally, the inequality $e^{mt}L_i\leq e^{-nt}X_{i+1}$ is equivalent to $t\leq q_i$ and the reverse inequality is equivalent to $t\geq q_i$. This proves the formula for $\lambda_1(t)$. Taking the logarithm yields \eqref{h1 formula}.
\end{proof}

A direct calculation gives
\begin{equation}\label{times diff}
    q_i-p_i=\frac{1}{d}\log \frac{X_{i+1}}{X_i}, \quad p_{i+1}-q_i=\frac{1}{d}\log \frac{L_i}{L_{i+1}}\leq D, \ D=\frac{\log M_2}{d}.
\end{equation}

We shall use the following standard exponential-growth estimate for best approximations; see \cite{BL}. Although the argument is stated for the supremum norm, it uses only reduction modulo 3 and the triangle inequality, and therefore applies to the Euclidean norms used here.

\begin{lemma}\label{growth lem}
One may take $S=3^d$ so that
\begin{equation}\label{growth eq}
    X_{i+S}\geq 2X_i
\end{equation}
for all sufficiently large $i$.
\end{lemma}

\begin{lemma}\label{long descent}
    For every $B>0$, there exists $T_B=T(B,d,M_2)>0$ such that, on a sufficiently late portion of the graph of $h_1$,
    \begin{equation}\label{long descent ineq}
        t-s\geq T_B \implies h_1(t)-h_1(s)> -n(t-s)+B.
    \end{equation}
    One may take
    \begin{equation}\label{T_B}
        N_B=S\left( \frac{B}{\log 2} +1 \right),\quad T_B=\frac{B}{d}+(N_B+3)D+1.
    \end{equation}
\end{lemma}
\begin{proof}
    We work on a sufficiently late portion of the graph such that \eqref{times diff}, \eqref{growth eq}, \eqref{h1 formula} hold. Take $\Delta=t-s\geq T_B$.

    Let $U$ be the sum of the lengths of intervals on $[s,t]$ on which $h_1$ is increasing, and let $U_d$ be the sum of the lengths of the intervals inside $[s,t]$ on which it is decreasing. Then $U+U_d=\Delta$. The slope is $m$ on each increasing piece and $-n$ on each decreasing piece. This gives
    \[
    h_1(t)-h_1(s)=mU-nU_d=dU-n\Delta.
    \]
    If \eqref{long descent ineq} fails, then $h_1(t)-h_1(s)\leq -n\Delta+B$. Therefore,
    \[
    U\leq \frac{B}{d}.
    \]

    Let $N$ be the number of complete increasing intervals $[p_i,q_i]$ inside $[s,t]$. Their indices are consecutive. Write them as $k,k+1,\dots,k+N-1$. Their total length is at most $U$. \eqref{times diff} yields
    \[
    dU\geq d\sum_{i=k}^{k+N-1}(q_i-p_i)=\sum_{i=k}^{k+N-1}\log \frac{X_{i+1}}{X_i}=\log \frac{X_{k+N}}{X_k}.
    \]
    Put $\nu=\lfloor N/S\rfloor$. Applying \eqref{growth eq} multiple times yields $X_{k+\nu S}\geq 2^\nu X_k$. Since $\nu S\leq N$ this gives
    \[
    X_{k+N}\geq X_{k+\nu S}\geq 2^\nu X_k.
    \]
    After taking the logarithm, we get
    \[
    \left\lfloor \frac{N}{S} \right\rfloor\log 2\leq \log \frac{X_{k+N}}{X_k}\leq dU\leq B.
    \]
    The inequality $N/S<\lfloor N/S\rfloor+1$ gives
    \[
    N<S\left( \frac{B}{\log 2} +1 \right)=N_B.
    \]

    There are at most $N+2$ increasing pieces. Since there is a decreasing piece between each one and, possibly, a decreasing piece at either end, there are no more than $N+3$ decreasing pieces. Each decreasing piece has length at most $D$ by \eqref{times diff}. Therefore
    \[
    U_d\leq (N+3)D\leq (N_B+3)D.
    \]
    Together with the estimate for $U$ this yields
    \[
    \Delta=U+U_d\leq \frac{B}{d}+(N_B+3)D=T_B-1<T_B.
    \]
    This contradicts $\Delta\geq T_B$.
\end{proof}

\subsection{Successive minima in a subspace}

Throughout this section let $V$ be a fixed rational subspace of dimension $r$. Define the lattice $\Gamma _V=V\cap \mathbb{Z}^d$. Let $\lambda_1(t),\dots,\lambda_r(t)$ denote the successive minima of the lattice $A_t\Gamma_V$ with respect to the norm $\|\cdot\|_\square$. In other words $\lambda_k(t)$ is the smallest number $\lambda$ such that there exist $k$ linearly independent vectors $z\in\Gamma_V$ satisfying $\|A_tz\|_\square \leq \lambda$. Let
\[
f_1(t)\leq \dots \leq f_r(t)
\]
denote the logarithms of the successive minima. Set
\[
S_k(t)=f_1(t)+\dots+f_k(t).
\]

\begin{lemma}\label{lem:fixed short subspace}
    The functions $f_k$ are continuous. For $t\geq s$,
    \begin{equation}\label{fk diff ineq}
        -n(t-s)\leq f_k(t)-f_k(s)\leq m(t-s).
    \end{equation}
    If $1\leq k<r$ and $f_k<f_{k+1}$ on an interval $J$, then the subspace
    \begin{equation}\label{Uk(t)}
        U_k(t)=\operatorname{span}_{\mathbb{R}}\{z\in\Gamma_V:\|A_tz\|_\square<e^{f_{k+1}(t)}\}
    \end{equation}
    has dimension $k$ and is independent of $t\in J$. The first $k$ minima of the lattice $A_t(U_k\cap \mathbb{Z}^d)$ equal $e^{f_1(t)},\dots,e^{f_k(t)}$.
\end{lemma}
\begin{proof}
    Put $\Delta=t-s\geq 0$. Then $A_t=D_\Delta A_s$. Thus, if $A_sz=(u,v)$, then
    \[
    \|A_tz\|_\square=\max(e^{-n\Delta}\|u\|_2,e^{m\Delta}\|v\|_2).
    \]
    This yields
    \[
    e^{-n\Delta}\|A_sz\|_\square\leq \|A_tz\|_\square\leq e^{m\Delta}\|A_sz\|_\square.
    \]
    We can take independent vectors $z_1,\dots,z_k\in\Gamma_V$ such that $\|A_sz_j\|_\square\leq \lambda_j(s)$ for $j=1,\dots,k$. Then by our estimate $\|A_tz_j\|_\square\leq e^{m\Delta}\lambda_k(s)$. Therefore
    \[
    \lambda_k(t)\leq e^{m\Delta}\lambda_k(s).
    \]
    Similarly we get the inequality
    \[
    e^{-n\Delta}\lambda_k(s)\leq \lambda_k(t).
    \]
    Taking logarithms in both inequalities gives \eqref{fk diff ineq}. Moreover, this inequality implies
    \[
    |f_k(t)-f_k(s)|\leq \max(m,n)|t-s|,
    \]
    which proves continuity.

    Now suppose that $f_k<f_{k+1}$ on some interval $J$, $1\leq k<r$. This implies that $\lambda_k(t)<\lambda_{k+1}(t)$ on $J$. Fix $t_0\in J$. Choose $k$ independent vectors $w_1,\dots,w_k\in\Gamma_V$ such that $\|A_{t_0}w_j\|_\square\leq \lambda_k(t_0)$ for $j=1,\dots,k$. Let $W=\operatorname{span}_{\mathbb{R}}\{w_1,\dots,w_k\}$. Since the functions $\|A_tw_j\|_\square$ are continuous, and $\lambda_k(t_0)<\lambda_{k+1}(t_0)$, there exists a neighborhood of $t_0$ such that for every $t$ from that neighborhood $\|A_tw_j\|_\square<\lambda_{k+1}(t)$. Therefore all these vectors $w_j$ lie in $U_k(t)$ by definition. The dimension of this space is $k$, since the maximal number of linearly independent vectors $z\in\Gamma_V$ satisfying $\|A_tz\|_\square<\lambda_{k+1}(t)$ is $k$. The dimension of $W$ is the same. Hence $W=U_k(t)$. We have proven that the function
    \[
    t \mapsto U_k(t)
    \]
    is locally constant on $J$. Since a locally constant function on a connected interval is globally constant, we have proven that $U_k(t)$ is independent of $t\in J$.

    Fix $1\leq j\leq k$. Let $\mu_j(t)$ be the $j\text{th}$ minimum of $A_t(U_k\cap \mathbb{Z}^d)$. Since all vectors from $U_k\cap \mathbb{Z}^d$ belong to $\Gamma_V$ we have
    \[
    \mu_j(t)\geq \lambda_j(t).
    \]
    For the reverse inequality choose $k$ independent vectors $w_1,\dots,w_k\in\Gamma_V$ attaining the first $k$ successive minima. For their norms we have
    \[
    \|A_tw_j\|_\square\leq \lambda_j(t)\leq \lambda_k(t)<\lambda_{k+1}(t).
    \]
    Therefore all of them belong to $U_k\cap \mathbb{Z}^d$, and
    \[
    \mu_j(t)\leq \lambda_j(t).
    \]
    These two inequalities give
    \[
    \mu_j(t)=\lambda_j(t),
    \]
    for $1\leq j\leq k$.
\end{proof}

Let $W\subset V$ be a rational $k-\text{dimensional}$ subspace of $V$. Choose a basis of the lattice $W\cap \mathbb{Z}^d$. Denote it by $w_1,\dots,w_k$. Put
\[
\Omega_W:=\Phi w_1\wedge \dots\wedge \Phi w_k.
\]
This vector is defined up to sign. We decompose it orthogonally as
\begin{equation}\label{orthogonal decomp}
    \Omega_W=\sum_{l=0}^k \Omega_{W,l}, \quad \Omega_{W,l}\in \bigwedge^{k-l}E\wedge \bigwedge^l F, \quad K_l(W)=\|\Omega_{W,l}\|_2.
\end{equation}
Since $\dim E=m$, $\dim F=n$ we have
\begin{equation}\label{dim bounds}
    K_l(W)\neq 0 \implies \max(0,k-m)\leq l\leq \min(n,k).
\end{equation}

\begin{lemma}
    \begin{equation}\label{Ht(W)}
        H_t(W):=\operatorname{covol}(A_t(W\cap \mathbb{Z}^d))=\left( \sum_{l=0}^k K_l(W)^2 e^{2(-nk+dl)t} \right)^{1/2}.
    \end{equation}
\end{lemma}
\begin{proof}
    We use the formula
    \[
    \operatorname{covol}(A_t(W\cap \mathbb{Z}^d))=\|A_tw_1\wedge\dots\wedge A_tw_k\|_2.
    \]
    The expression on the right equals
    \[
    \left\|\left(\bigwedge^k D_t\right)\Omega_W\right\|_2.
    \]
    By definition of $D_t$ it multiplies every vector from $E$ by $e^{-nt}$ and it multiplies every vector from $F$ by $e^{mt}$. Therefore
    \[
    \left(\bigwedge^k D_t\right) \Omega_{W,l}=e^{-nt(k-l)}e^{mtl}\Omega_{W,l}=e^{(-nk+dl)t}\Omega_{W,l}.
    \]
    The orthogonality of the decomposition finishes the proof.
\end{proof}

We shall also use
\begin{equation}\label{Ht wave}
    \tilde{H}_t(W):=\max_{l:K_l(W)>0} K_l(W)e^{(-nk+dl)t}, \quad \tilde{H}_t(W)\leq H_t(W)\leq \sqrt{k+1}\tilde{H}_t(W).
\end{equation}
In particular, $\log \tilde{H}_t(W)$ is a continuous and locally piecewise affine function.

\begin{lemma}\label{lem:comp of heights and succ minima}
    There exists a $\delta=\delta(d)\geq 1$ with the following property. Let $f_1^W(t),\dots,f_k^W(t)$ be the logarithms of the successive minima of $A_t(W\cap \mathbb{Z}^d)$ with respect to $\| \cdot \|_\square$. Then
    \begin{equation}\label{eq:comp of heights and succ minima}
        \left| \log \tilde{H}_t(W)-\sum_{j=1}^k f_j^W(t) \right|\leq \delta.
    \end{equation}
    This inequality holds for all $W,k,t$.
\end{lemma}
\begin{proof}
    Fix $W$ and $t$. Write
    \[
    U=A_tW, \quad \Lambda=A_t(W\cap \mathbb{Z}^d), \quad \lambda_j=e^{f_j^W(t)}, \quad \Pi = \prod_{j=1}^k \lambda_j.
    \]

    Define the unit ball
    \[
    B_\square=\{(u,v):\|u\|_2\leq 1,\|v\|_2\leq 1\}.
    \]
    Let $K=B_\square\cap U$. Then for $\lambda>0$, $z\in U$
    \[
    z\in \lambda K \iff \|z\|_\square \leq\lambda.
    \]
    Consequently, the successive minima of $K$ with respect to $\Lambda$ are $\lambda_1,\dots,\lambda_k$. The set $K$ is compact, convex, and symmetric about the origin.

    Let
    \[
    B_U(\rho)=\{z\in U:\|z\|_2\leq \rho\}.
    \]
    The inequalities \eqref{block and eucl norm comp} imply
    \[
    B_U(1)\subset K\subset B_U(\sqrt{2}).
    \]
    In particular, the $k$-dimensional volume of $K$ is positive and finite. Let $v_k$ be the volume of the Euclidean unit ball in $\mathbb{R}^k$, and put $V_K=\operatorname{vol}_k(K)$. The inclusion implies
    \begin{equation}\label{eq:vol comp}
    v_k\leq V_K\leq (\sqrt{2})^kv_k=2^{k/2}v_k.
    \end{equation}

    Applying Minkowski's second theorem to $K$ and the lattice $\Lambda$ yields
    \[
    \frac{2^k}{k!}\det \Lambda\leq V_K\Pi \leq 2^k \det \Lambda.
    \]
    By definition $\det \Lambda=H_t(W)$. Hence
    \[
    \frac{2^k}{k!V_K}\leq \frac{\Pi}{H_t(W)}\leq \frac{2^k}{V_K}.
    \]
    Using \eqref{eq:vol comp} yields
    \[
    \frac{2^{k/2}}{k!v_k}\leq \frac{\Pi}{H_t(W)}\leq \frac{2^k}{v_k}.
    \]
    \eqref{Ht wave} in the form
    \[
    1\leq \frac{H_t(W)}{\tilde{H}_t(W)}\leq\sqrt{k+1}
    \]
    and
    \[
    \frac{\Pi}{\tilde{H}_t(W)}=\frac{\Pi}{H_t(W)}\frac{H_t(W)}{\tilde{H}_t(W)}
    \]
    give
    \[
    a_k\leq \frac{\prod_{j=1}^k \lambda_j}{\tilde{H}_t(W)}\leq b_k, \quad a_k=\frac{2^{k/2}}{k!v_k},\quad b_k=\frac{2^k\sqrt{k+1}}{v_k}.
    \]
    After taking the logarithm we obtain
    \[
    \log a_k\leq \sum_{j=1}^k f_j^W(t)-\log \tilde{H}_t(W)\leq \log b_k.
    \]
    Choose $\delta =\max\limits_{1\leq k\leq d} \{1,|\log a_k|,|\log b_k|\}$. Then \eqref{eq:comp of heights and succ minima} holds.
\end{proof}

\begin{lemma}\label{lem:min height of rank k}
    Define
    \begin{equation}\label{eq:Gk(t)}
        G_k(t)=\min_{\substack{W\subseteq V\ \text{rational}\\ \dim W=k}}\log \tilde{H}_t(W).
    \end{equation}
    This minimum is attained, and the functions $G_k$ are continuous and locally piecewise affine. Moreover,
    \begin{equation}\label{eq:Gk Sk comp}
        |G_k(t)-S_k(t)|\leq \delta, \quad G_1(t)=f_1(t).
    \end{equation}
    If $k<r$ and $f_{k+1}(t)-f_k(t)>2\delta$, then every subspace attaining the minimum in \eqref{eq:Gk(t)} is equal to $U_k(t)$.
\end{lemma}
\begin{proof}
    Fix a compact time interval $J$. Choose any  rational $k-\text{dimensional}$ subspace $W_0\subset V$. Put
    \[
    B_J=\max_{t\in J} \tilde{H}_t(W_0)>0, \quad R_J=\max_{t\in J}\left\| \left(\bigwedge^k A_t \right)^{-1}\right\|_2.
    \]
    Here $B_J$ is finite since the function $t\mapsto\tilde H_t(W_0)$ is continuous, and $R_J$ is finite since all the matrix entries of $A_t^{-1}$ continuously depend on $t$.

    Only a plane $W$ with $\tilde{H}_{t_0}(W)\leq B_J$ at some $t_0\in J$ can attain the minimum anywhere on $J$. Choose a basis of the lattice $W\cap \mathbb{Z}^d$. Denote this basis by $w_1,\dots,w_k$. Let $\eta_W=w_1\wedge\dots\wedge w_k$. Then at time $t_0$ we have
    \[
    \|\eta_W\|_2\leq R_J \left\| \left(\bigwedge^k A_{t_0} \right)\eta_W \right\|_2=R_JH_{t_0}(W)\leq R_J \sqrt{k+1}\tilde{H}_{t_0}(W)\leq R_J\sqrt{k+1}B_J.
    \]
    This inequality bounds every coordinate of $\eta_W$. Since its coordinates are integers, it can take only finitely many values. Furthermore, each such vector defines its plane
    \[
    W=\{z\in\mathbb{R}^d:z\wedge\eta_W=0\}.
    \]
    Therefore, the number of planes that can compete with $W_0$ anywhere on $J$ is finite. Consequently, the minimum in \eqref{eq:Gk(t)} can be replaced by a minimum over a finite family. Hence it is attained.

    As mentioned before, each function $\log \tilde{H}_t(W)$ is continuous and locally piecewise affine. So is the minimum of a finite number of them.

    Now we prove the inequality in \eqref{eq:Gk Sk comp}. Fix an arbitrary $t$ and a $k-\text{dimensional}$ rational subspace $W\subseteq V$. Write $\lambda_j^W=e^{f_j^W(t)}$. Restricting $\Gamma_V$ to $W\cap\mathbb{Z}^d$ leaves fewer vectors. Hence $\lambda_j\leq \lambda_j^W$ and $f_j(t)\leq f_j^W(t)$. \Cref{lem:comp of heights and succ minima} gives
    \[
    \log \tilde{H}_t(W)\geq \sum_{j=1}^k f_j^W(t)-\delta\geq \sum_{j=1}^k f_j(t)-\delta=S_k(t)-\delta.
    \]
    Taking the minimum over all such $W$ proves one side of the inequality. Now we show the reverse inequality. Consider
    \[
    E_j=\{z\in\Gamma_V:\|A_tz\|_\square\leq \lambda_j\}.
    \]
    By definition of $\lambda_j$ we have $\dim \operatorname{span}_{\mathbb{R}} E_j\geq j$. Now we construct a rational $k-\text{dimensional}$ subspace $W\subset V$. First, choose nonzero $z_1\in E_1$. Suppose we have already chosen independent $z_1,\dots,z_{j-1}$, where $2\leq j\leq k$ for which $\|A_tz_i\|_\square\leq \lambda_i$. Put
    \[
    P_{j-1}=\operatorname{span}_{\mathbb{R}}\{z_1,\dots,z_{j-1}\}.
    \]
    Its dimension is $j-1$. Therefore it does not contain $E_j$. Choose $z_j\in E_j\setminus P_{j-1}$. Then the vectors $z_1,\dots,z_j$ are independent and $\|A_tz_j\|_\square\leq \lambda_j$. This way we can construct a system of independent vectors $z_1,\dots,z_k$ for which
    \[
    \|A_tz_i\|_\square\leq\lambda_i, \quad (1\leq i\leq k).
    \]
    Put
    \[
    W=\operatorname{span}_{\mathbb{R}}\{z_1,\dots,z_k\}.
    \]
    By construction we have $\lambda_j^W\leq \lambda_j$. The reverse inequality is obvious. Hence $\lambda_j^W= \lambda_j$, and $f_j^W=f_j$. Now \cref{lem:comp of heights and succ minima} gives
    \[
    G_k(t)\leq \log \tilde{H}_t(W)\leq \sum_{j=1}^k f_j^W(t)+\delta =\sum_{j=1}^k f_j(t)+\delta=S_k(t)+\delta.
    \]
    This concludes the proof of the inequality in \eqref{eq:Gk Sk comp}.

    To show that $G_1(t)=f_1(t)$ we fix a line $W$ generated by a primitive vector $w\in \Gamma_V$. Write $\Phi w=(u,v)$. Formula \eqref{Ht wave} shows that in the case $k=1$ we have
    \[
    \tilde{H}_t(W)=\max(e^{-nt}\|u\|_2,e^{mt}\|v\|_2)=\|A_tw\|_\square.
    \]
 Therefore, minimizing this expression is the same as minimizing $\|A_tw\|_\square$ with $0\neq w\in\Gamma_V$. Hence $G_1(t)=f_1(t)$.

    Now suppose $k<r$ and $f_{k+1}(t)-f_k(t)>2\delta$. \Cref{lem:fixed short subspace} provides a $k-\text{dimensional}$ subspace $U_k(t)$. If $W\neq U_k(t)$, then $\lambda_k^W\geq\lambda_{k+1}$. Otherwise, $W$ would contain $k$ independent integer vectors that have norms after applying $A_t$ that are strictly less than $\lambda_{k+1}$. Therefore, all of them belong to $U_k(t)$ by definition. The equality of their dimensions gives $W=U_k(t)$, which is a contradiction. Taking logarithms gives $f_k^W(t)\geq f_{k+1}(t)$. For $j<k$ we have $f_j^W(t)\geq f_j(t)$. This, \cref{lem:comp of heights and succ minima} and \eqref{eq:Gk Sk comp} yield
    \[
    \log \tilde{H}_t(W)\geq \sum_{j=1}^k f_j^W(t)-\delta \geq S_k(t)+f_{k+1}(t)-f_k(t)-\delta>S_k(t)+\delta\geq G_k(t).
    \]
    This proves that no subspace $W\neq U_k(t)$ can attain the minimum.
\end{proof}

\subsection{The potential and its properties}\label{sec:potential}

Choose $t_*$ such that \cref{long descent} holds with $B=2\delta$ whenever $t>s\geq t_*$ and $t-s\geq T_B$. Let $I=[u_0,u_1],u_0\geq t_*$ be an interval, with the right end possibly being infinite, on which $f_1(t)=h_1(t)$. For the application to \cref{th:block lemma} we take $I=[q_a,q_{b-1}]$. Here is why the equality holds on this interval. Since $\Gamma_V\subset\mathbb{Z}^d$ we have $f_1\geq h_1$. On the other hand, the proof of \eqref{h1 formula} shows that at every $t$ on this interval the minimum of $\|A_tz\|_\square$ is attained by one of the vectors $z_a,\dots,z_b$. All these vectors lie in $V$ by hypothesis in \cref{th:block lemma}. This gives $f_1\leq h_1$. Thus, they are equal on $I$.

Fix
\begin{equation}\label{eq:T,c}
    T=T_{2\delta},\quad c=3r\delta+dT+1,
\end{equation}
and define
\begin{equation}\label{eq:potential}
    P(t):=\min_{1\leq k\leq r} \left\{ \frac{G_k(t)}{k}+(r-k)c \right\}.
\end{equation}
A rank $k$ is called active if it attains this minimum. By \cref{lem:min height of rank k} the potential is continuous and locally piecewise affine, since it is the minimum of a finite number of such functions.

\begin{proposition}
    The following inequalities hold
    \begin{equation}\label{eq:potential comp}
        f_1(t)-\delta\leq P(t)\leq f_1(t)+(r-1)c.
    \end{equation}
\end{proposition}
\begin{proof}
    Since $f_j(t)\geq f_1(t)$ for $j\leq k$, we have $S_k(t)\geq kf_1(t)$. By \cref{lem:min height of rank k} for every rank $k$ we have
    \[
    \frac{G_k(t)}{k}+(r-k)c\geq \frac{S_k(t)}{k}-\frac{\delta}{k}+(r-k)c\geq f_1(t)-\delta.
    \]
    Taking the minimum proves the left side of \eqref{eq:potential comp}. As for the right side, the branch of rank one is admissible. Since $G_1(t)=f_1(t)$, we have
    \[
    P(t)\leq f_1(t)+(r-1)c.
    \]
\end{proof}

\begin{lemma}\label{lem:active rank and gap}
    If a rank $k<r$ is active at a time $t$, then
    \begin{equation}\label{eq:active rank gap ineq}
        f_{k+1}(t)-f_k(t)\geq 2c-4k\delta>dT+2\delta.
    \end{equation}
\end{lemma}
\begin{proof}
    Fix the time $t$. Let
    \[
    a_j=\frac{S_j}{j}, \quad \varepsilon_j=G_j-S_j\implies |\varepsilon_j|\leq\delta.
    \]
    First, assume that $2\leq k<r$. Compare the active branch to two neighboring branches. Comparing with $k-1$ gives
    \[
    \frac{G_k}{k}+(r-k)c\leq \frac{G_{k-1}}{k-1}+(r-k+1)c\implies \frac{G_k}{k}-\frac{G_{k-1}}{k-1}\leq c.
    \]
    Comparing with $k+1$ gives
    \[
    \frac{G_k}{k}+(r-k)c\leq\frac{G_{k+1}}{k+1}+(r-k-1)c\implies \frac{G_{k+1}}{k+1}-\frac{G_k}{k}\geq c.
    \]
    Substituting $G_j/j=a_j+\varepsilon_j/j$ and using $|\varepsilon_j/j|\leq\delta$ yields
    \begin{equation}\label{eq:ak comp}
    a_k-a_{k-1}\leq c+2\delta, \quad a_{k+1}-a_k\geq c-2\delta.
    \end{equation}
    One can easily check
    \[
    f_k-a_k=(k-1)(a_k-a_{k-1}),\quad f_{k+1}-a_k=(k+1)(a_{k+1}-a_k).
    \]
    Subtracting the two expressions and using \eqref{eq:ak comp} gives
    \[
    f_{k+1}-f_k=(k+1)(a_{k+1}-a_k)-(k-1)(a_k-a_{k-1})\geq (k+1)(c-2\delta)-(k-1)(c+2\delta)=2c-4k\delta.
    \]
    In the case $k=1$ we only compare this branch with the branch of rank $2$. Similarly this comparison gives $a_2-a_1\geq c-2\delta$. Since $a_1=f_1$, and $2a_2=f_1+f_2$, we get $f_2-f_1=2(a_2-a_1)\geq 2c-4\delta$, as required. The comparison $2c-4k\delta>dT+2\delta$ can be checked easily.
\end{proof}

\begin{lemma}\label{lem:l=0}
    Let $t\in I$ and $t\geq u_0+T$. Suppose a rank $k$ is active and $W$ attains the minimum in $G_k(t)$. Then the term with $l=0$ does not attain the maximum in \eqref{Ht wave}, even in a tie with other terms.
\end{lemma}
\begin{proof}
    The conditions $t\geq u_0+T$ and $t\in I$ imply $[t-T,t]\subseteq I$. First suppose $k<r$. Let $g=f_{k+1}-f_k$. For $v_2\geq v_1$ the inequalities in \eqref{fk diff ineq} give
    \[
    -d(v_2-v_1)\leq g(v_2)-g(v_1)\leq d(v_2-v_1).
    \]
    By lemma \ref{lem:active rank and gap}, $g(t)>dT+2\delta$. Hence for every $v\in [t-T,t]$,
    \[
    g(v)\geq g(t)-d(t-v)>2\delta>0.
    \]
    Lemma \ref{lem:fixed short subspace} says that $U_k(v)$ is a fixed subspace for all $v$ in this interval, and that the first $k$ minima of $A_v(U_k(v)\cap\mathbb{Z}^d)$ are the same as those of $V$. For every $v$ in this interval, Lemma \ref{lem:min height of rank k} says that the minimizing subspace $W=U_k(v)$ is equal to $U_k(t)$. Applying Lemma \ref{lem:comp of heights and succ minima} to this $W$ at every time $v$ gives
    \begin{equation}\label{eq:lem Sk H wave diff}
    |S_k(v)-\log \tilde{H}_v(W)|\leq\delta.
    \end{equation}
    If $k=r$ this inequality holds as well, since $W=V$.

    Now we assume that the term with $l=0$ does attain the maximum at time $t$. Then $K_0(W)>0$. For every $l>0$ consider the following expression
    \[
    \frac{K_l(W)}{K_0(W)}e^{dlv}.
    \]
    It increases on $[t-T,t]$ and is at most $1$ by our assumption. Hence $l=0$ also attains the maximum for all $v$ from this interval. Therefore,
    \[
    \log \tilde{H}_v(W)=\log K_0(W)-nkv.
    \]
    By \eqref{eq:lem Sk H wave diff},
    \[
    S_k(t)\leq \log \tilde{H}_t(W)+\delta,
    \]
    and
    \[
    S_k(t-T)\geq\log\tilde{H}_{t-T}(W)-\delta.
    \]
    Subtracting them yields
    \[
    S_k(t)-S_k(t-T)\leq -nkT+2\delta.
    \]
    For each $j=2,\dots, k$ the inequality \eqref{fk diff ineq} gives $f_j(t)-f_j(t-T)\geq -nT$. Therefore
    \[
    f_1(t)-f_1(t-T)=S_k(t)-S_k(t-T)-\sum_{j=2}^k (f_j(t)-f_j(t-T))\leq -nkT+2\delta +n(k-1)T=-nT+2\delta.
    \]
    For $k=1$ the sum is empty, and the same estimate holds. Since $f_1=h_1$ on $I$, By  \cref{long descent}
    \[
    h_1(t)-h_1(t-T)>-nT+2\delta.
    \]
    This is a contradiction.
\end{proof}

\begin{lemma}\label{lem:lower bound for slope}
    Set $\mu=\min(r-1,m)$ and
    \[
    \sigma=-n+\frac{d}{\mu+1}.
    \]
    For $s\leq t$ in $I$, we have
    \begin{equation}\label{eq:lower bounde for slope}
        f_1(t)-f_1(s)\geq \sigma(t-s)-C_0,\quad C_0=(r-1)c+\delta+dT.
    \end{equation}
    Notice that $C_0$ depends only on $m,n,M_2,r$.
\end{lemma}
\begin{proof}
    As mentioned before, $P$ is a continuous locally piecewise affine function. By its definition, at each point $v$
    \[
    P(v)=\frac{\log K_l(W)+(-nk+dl)v}{k}+(r-k)c.
    \]
    Therefore, on each interval where $P$ is an affine function, its slope is
    \[
    -n+\frac{dl}{k}.
    \]

    Assume that $u_0+T\leq v\leq u_1$. In this case, by \cref{lem:l=0} we have $l\neq 0$ in the formula for $P(v)$. Moreover, \eqref{dim bounds} implies that $l\geq k-m$. Hence,
    \[
    l\geq\max(1,k-m).
    \]
    We claim
    \[
    \frac{l}{k}\geq \frac{1}{\min(r,m+1)}=\frac{1}{\mu+1}.
    \]
    To see this, let $q=\min(r,m+1)=\mu+1$. If $k\leq q$, then
    \[
    \frac{l}{k}\geq \frac{1}{k}\geq\frac{1}{q}.
    \]
    If $k>q$, then $q=r$ is impossible, since $k\leq r$. Hence, $q=m+1$ and $k>m+1$. In this case
    \[
    \frac{l}{k}\geq\frac{k-m}{k}, \quad \frac{k-m}{k}-\frac{1}{m+1}=\frac{m(k-m-1)}{k(m+1)}>0.
    \]
    This proves the inequality in both cases. Therefore, in this case the following estimate  holds for the slope
    \[
    -n+\frac{dl}{k}\geq\sigma.
    \]

    Now we consider different cases. Put $v_0=u_0+T$.
    \begin{enumerate}
        \item $v_0\leq s\leq t$.
        
        In this case, divide $[s,t]$ into finitely many pieces $[v_i,v_{i+1}]$, such that $P$ is affine on each subinterval. On each such piece we have
        \[
        P(v_{i+1})-P(v_i)\geq\sigma (v_{i+1}-v_i).
        \]
        Summing this yields
        \[
        P(t)-P(s)\geq\sigma (t-s).
        \]
        \eqref{eq:potential comp} gives
        \[
        f_1(t)\geq P(t)-(r-1)c,\quad f_1(s)\leq P(s)+\delta.
        \]
        Hence,
        \[
        f_1(t)-f_1(s)\geq\sigma (t-s)-A,\quad A=(r-1)c+\delta.
        \]

        \item $s<v_0<t$.

        We apply \eqref{fk diff ineq} on $[s,v_0]$ and the previous case on $[v_0,t]$. This gives
        \[
        f_1(t)-f_1(s)\geq -n(v_0-s)+\sigma (t-v_0)-A=\sigma (t-s)-A-(n+\sigma)(v_0-s).
        \]
        Since $v_0-s\leq T$ and $n+\sigma=d/(\mu+1)\leq d$, the last loss is at most $dT$.

        \item $s<t\leq v_0$.

        In this case $t-s\leq T$. We apply \eqref{fk diff ineq}
        \[
        f_1(t)-f_1(s)\geq -n(t-s)=\sigma (t-s)-(n+\sigma)(t-s)\geq \sigma (t-s)-dT.
        \]
        The constant $C_0=A+dT$ covers all cases.
    \end{enumerate}
\end{proof}

\begin{lemma}\label{lem:short ind sys}
    Suppose rank $k$ is active for the potential at time $t$. Then
    \begin{equation}\label{eq:fk-f1 ineq}
        f_k(t)-f_1(t)\leq k(k-1)c+\delta.
    \end{equation}
    Set
    \[
    C_k=e^{k(k-1)c+\delta}.
    \]
    There exist $k$ independent vectors $v_1,\dots,v_k\in \Gamma_V$ such that
    \begin{equation}\label{eq:vj comp}
    \|A_tv_j\|_\square\leq C_k e^{f_1(t)}
    \end{equation}
    If $k<r$, their span coincides with the minimizing subspace $W$.
\end{lemma}
\begin{proof}
    We compare the active rank $k$ to rank $1$,
    \[
    \frac{G_k(t)}{k}+(r-k)c\leq f_1(t)+(r-1)c.
    \]
    This gives
    \[
    G_k(t)\leq kf_1(t)+k(k-1)c.
    \]
    The comparison in \eqref{eq:Gk Sk comp} gives $S_k(t)\leq G_k(t)+\delta$. Thus, $S_k(t)\leq kf_1(t)+k(k-1)c+\delta$. Combining it with $S_k(t)\geq (k-1)f_1(t)+f_k(t)$ and rearranging proves \eqref{eq:fk-f1 ineq}.

    We have $k$ independent vectors $v_1,\dots,v_k\in\Gamma_V$, for which $\|A_tv_j\|_\square\leq e^{f_k(t)}$. Exponentiating \eqref{eq:fk-f1 ineq} gives
    \[
    \|A_tv_j\|_\square\leq C_ke^{f_1(t)},
    \]
    which is exactly \eqref{eq:vj comp}.

    Now we identify the span of these vectors in the case $k<r$. \Cref{lem:active rank and gap} gives $f_{k+1}(t)-f_k(t)> 2\delta$. \Cref{lem:min height of rank k} states that the minimizing subspace equals $U_k(t)$, which is exactly the span of the chosen vectors.
\end{proof}

\section{\texorpdfstring
  {Proof of \Cref{th:block lemma}}
  {Proof of Theorem \ref*{th:block lemma}}}

\begin{proof}[Proof of \Cref{th:block lemma}]
We retain the notation from the previous section. Take $a$ large enough so that $q_a\geq t_*$, and set $I=[q_a,q_{b-1}]$. As mentioned before this $I$ is admissible. Note that $r=\dim V\geq 2$, since consecutive best approximation vectors are independent.

Formula \eqref{h1 formula} gives
\[
f_1(q_i)=\log L_i +mq_i.
\]
Also,
\[
m-\sigma=m+n-\frac{d}{\mu+1}=\frac{d\mu}{\mu+1}.
\]
Using the definition of $q_i$ in \eqref{pi qi}, we obtain
\[
f_1(q_i)-\sigma q_i=\log L_i+(m-\sigma)q_i=\log L_i+\frac{\mu}{\mu+1}(\log X_{i+1}-\log L_i)=\frac{\log L_i+\mu\log X_{i+1}}{\mu+1}.
\]
Hence
\begin{equation}\label{eq:f-sigma}
    f_1(q_i)-\sigma q_i=\frac{1}{\mu+1}\log(L_iX_{i+1}^\mu).
\end{equation}

Apply \cref{lem:lower bound for slope} with $s=q_i$, $t=q_j$:
\[
(f_1(q_j)-\sigma q_j)-(f_1(q_i)-\sigma q_i)\geq -C_0.
\]
Now substitute \eqref{eq:f-sigma} into this inequality:
\[
\log(L_iX_{i+1}^\mu)\leq\log(L_jX_{j+1}^\mu)+(\mu+1)C_0.
\]
After exponentiation, this becomes \eqref{block lemma ineq} with $C_*^{(r)}=e^{(\mu+1)C_0}$. To make this constant independent of $r$, take $C_*=\max_{2\leq r\leq d}C_*^{(r)}$. 
\end{proof}

\section{\texorpdfstring
  {Proof of \Cref{problem 1 sol}}
  {Proof of Theorem \ref*{problem 1 sol}}}

\begin{proof}[Proof of \Cref{problem 1 sol}]
    In \cref{th:block lemma} take $V=\mathbb{R}^d$. Then $r=d=m+n$ and $r-1\geq m$. Hence $\mu=\min(m,r-1)=m$. Fix $i_0$ sufficiently large so that the block lemma applies. Choose an index $i\geq i_0+1$. \eqref{block lemma ineq} gives
    \[
    L_{i-1}X_i^m\geq C_*^{-1}L_{i_0}X_{i_0+1}^m:=B_0>0.
    \]
    Condition \eqref{cond C} implies $L_i\geq M_2^{-1}L_{i-1}$. Hence
    \[
    L_iX_i^m\geq B_0M_2^{-1}:=c_\xi.
    \]
    Taking the logarithm and rearranging gives
    \[
    \frac{-\log L_i}{\log X_i}\leq m-\frac{\log c_\xi}{\log X_i}.
    \]
    Taking the upper limit and using \eqref{eq:exp formulas} yields
    \[
    \omega(\xi)\leq m,
    \]
    as required.
\end{proof}

The following corollary is needed to prove \cref{problem 3 sol}. Recall the notation introduced before it.

\begin{corollary}\label{C=BC ord exp}
    For $n\geq 2$
    \begin{equation}\label{sigma omega C}
    \Sigma_\omega(\mathcal{C}_{m,n})=\left[ \frac{m}{n},m \right].
    \end{equation}
    For $m,n\geq 2$
    \begin{equation}\label{sigma omega BC}
        \Sigma_\omega(\mathcal{BC}_{m,n})=\left[ \frac{m}{n},m \right].
    \end{equation}
\end{corollary}
\begin{proof}
    Let $n\geq 2$. Assume $\xi$ is totally irrational and \eqref{cond C} holds. It is well known that $\omega(\xi)\geq \frac{m}{n}$. \Cref{problem 1 sol} proves that under condition \eqref{cond C} $\omega(\xi)\leq m$. Hence the whole spectrum is contained in $\left[ \frac{m}{n},m \right]$. In \cite[Theorem~2.6(2)]{MMS} it is proven that every exponent from this interval is admissible. Thus \eqref{sigma omega C} is proven.

    Let $m,n\geq 2$. The proof that no exponent outside the interval $\left[ \frac{m}{n},m \right]$ is admissible is the same. Once again, \cite[Theorem~2.6(3)]{MMS} proves that every such exponent is admissible. This proves \eqref{sigma omega BC}.
\end{proof}

\section{\texorpdfstring
  {Proof of \Cref{tail theorem}}
  {Proof of Theorem \ref*{tail theorem}}}

We first prove an auxiliary lemma.

\begin{lemma}\label{lem:largest err deg}
    Let $Q$ be a rational subspace of dimension $r$, let $U=\Phi Q$, and $\rho=\dim (U\cap E)$. Then the largest index $l$ for which $K_l(Q)>0$ is $r-\rho$.
\end{lemma}
\begin{proof}
    Let $\pi_F:E\oplus F\to F$ be the projection operator. Restricted to $U$ its kernel is $U\cap E$ of dimension $\rho$. The rank-nullity theorem gives
    \[
    q:=\dim \pi_F(U)=\dim U-\dim (U\cap E)=r-\rho.
    \]
    All $F-\text{factors}$ in an expansion of an exterior product of vectors from $U$ belong to the space $\pi_F(U)$ of dimension $q$. If a term has more than $q$ such factors it vanishes, since these vectors are dependent. Therefore, every component of degree $l>q$ vanishes.

    Now we prove that degree $q$ actually occurs. To show this, choose a basis $u_1,\dots,u_\rho$ of $U\cap E$. Choose $v_1,\dots, v_q\in U$ such that $\pi_Fv_1,\dots,\pi_Fv_q$ form a basis of $\pi_F(U)$. Together $u_i$ and $v_j$ form a basis of $U$. Let $v_j=\pi_Ev_j+\pi_Fv_j$. The $u_i$ contribute no $F-\text{factors}$. Consider the wedge product
    \[
    u_1\wedge\dots\wedge u_\rho\wedge v_1\wedge\dots\wedge v_q.
    \]
    There is only one term in this product with exactly $q$ $F-\text{factors}$. That is
    \[
    u_1\wedge\dots\wedge u_\rho\wedge \pi_F v_1\wedge\dots\wedge \pi_F v_q.
    \]
    This wedge product is nonzero because the $u_i$ are linearly independent in $E$, the vectors $\pi_Fv_j$ are linearly independent in $F$, and $E\cap F=\{0\}$.

    In the definition of $\Omega_Q$ we might use a different basis of $U$ to find the wedge product. However, the product will only be multiplied by a nonzero number, that is the determinant of the change-of-basis matrix. Therefore, the component of degree $q$ stays nonzero.
\end{proof}

\begin{proof}[Proof of \Cref{tail theorem}]
    First assume $r=\dim Q\leq m$. Since consecutive best approximation vectors are independent, $r\geq 2$. Similarly to the proof of \cref{problem 1 sol} we get
    \[
    L_iX_i^{r-1}\geq c,
    \]
    where $c>0$ depends on $\xi,Q,M_2$. Since $X_i\geq 1$, and $m\geq r$,
    \[
    L_iX_i^{m-1}=(L_iX_i^{r-1})X_i^{m-r}\geq c.
    \]

    Now assume $r=m+1$. Since $n\geq 2$, we have $r=m+1<d$, so $Q$ is a proper rational subspace. Condition \eqref{cond2} says $\mathcal{L}_\xi\not\subset Q$. $\mathcal{L}_\xi$ has dimension $m$, and $\Phi\mathcal{L}_\xi=E$. Invertibility of $\Phi$ gives $\Phi Q\cap E=\Phi(Q\cap \mathcal{L}_\xi)$. Hence
    \[
    \rho=\dim(\Phi Q\cap E)=\dim(Q\cap \mathcal{L}_\xi)\leq m-1.
    \]
    \Cref{lem:largest err deg} says that the largest nonzero error degree is $l_*=m+1-\rho$. It satisfies $l_*\geq 2$.

    For any smaller degree with $K_l(Q)>0$ the ratio of its term in \eqref{Ht wave} to the $\text{degree-}l_*$ term is
    \[
    \frac{K_l(Q)}{K_{l_*}(Q)}e^{-d(l_*-l)t}\to 0, \quad (t\to\infty).
    \]
    After some time $t_Q$ all these ratios are less than $1$ simultaneously. Thus
    \[
    \log\tilde{H}_t(Q)=\log K_{l_*}(Q)+(-nr+dl_*)t.
    \]

    We use the potential introduced in \cref{sec:potential} with $V=Q$. $f_1=h_1$ on a sufficiently late half-line $[u_0,+\infty)$, since all sufficiently late best approximation vectors lie in $Q$. Choose $t_1\geq u_0+T, t_1\geq t_Q$. Let $t\geq t_1$. On an affine piece of $P(t)$, if the active rank is $k\leq m$, \cref{lem:l=0} gives $l\geq 1$, so its slope is at least
    \[
    -n+\frac{d}{k}\geq -n+\frac{d}{m}.
    \]
    If the rank $k=m+1$ is active, then the only subspace is $Q$ itself. This gives the same lower bound for the slope
    \[
    -n+\frac{dl_*}{m+1}\geq-n+\frac{2d}{m+1}\geq-n+\frac{d}{m}.
    \]
    Hence the slope of $P(t)$ on every sufficiently late affine piece is at least
    \[
    \sigma_Q=-n+\frac{d}{m}.
    \]
    This gives
    \[
    P(t)\geq P(t_1)+\sigma_Q(t-t_1).
    \]
    Since $f_1(t)\geq P(t)-(r-1)c$, we obtain $f_1(t)-\sigma_Qt\geq P(t_1)-\sigma_Q t_1-(r-1)c$. Let $C_Q=\max(0,-P(t_1)+\sigma_Q t_1+(r-1)c)\geq 0$. Therefore,
    \[
    f_1(t)-\sigma_Q t\geq -C_Q.
    \]
    The identity \eqref{eq:f-sigma} holds for any $\mu$ and corresponding $\sigma=-n+d/(\mu+1)$. Take $\mu=m-1$, for which $\sigma=\sigma_Q$. For all sufficiently large $i$ it gives
    \[
    \log (L_iX_{i+1}^{m-1})=m(f_1(q_i)-\sigma_Qq_i)\geq-mC_Q\implies L_iX_{i+1}^{m-1}\geq e^{-mC_Q}>0.
    \]
    Condition \eqref{cond C} gives
    \[
    L_iX_i^{m-1}\geq M_2^{-1}e^{-mC_Q}=C.
    \]
    Taking the logarithm and using \eqref{eq:exp formulas} gives
    \[
    \omega(\xi)\leq m-1.
    \]
\end{proof}

\section{\texorpdfstring
  {Proof of \Cref{problem 2 sol}}
  {Proof of Theorem \ref*{problem 2 sol}}}

Throughout this section, $m,n\geq 2$, $V=\mathbb{R}^d$, and conditions \eqref{cond1}, \eqref{cond2}, \eqref{cond C} are assumed. Let $h(t)=h_1(t)=f_1(t)$ for $t\geq u_0\geq t_*$. Constants denoted by $C$ may change from line to line and are allowed to depend on $\xi,m,n,M_2$.

We introduce the following constants
\begin{equation}\label{eq:constants}
    \begin{split}
        \Delta_m=m^2+m+1,\quad s=-n+\frac{d}{m+1}=\frac{m(1-n)}{m+1},\\ u=-n+\frac{d}{m},\quad \kappa=m(2-n),\quad \beta_*=\frac{m(m+1-mn)}{\Delta_m}.
    \end{split}
\end{equation}
It can be easily checked that these constants satisfy
\begin{equation}\label{eq:const relations}
    s<\beta_*<u,\quad \frac{\kappa}{m+2}-\beta_*=\frac{m(m-1)d}{(m+2)\Delta_m}>0,\quad u-\beta_*=\frac{d}{m\Delta_m}>0.
\end{equation}
Moreover, the desired bound for the uniform exponent corresponds to $\beta_*$ in the following way
\begin{equation}\label{eq:beta*}
    \frac{m-\beta_*}{n+\beta_*}=\frac{m^2}{m+1}.
\end{equation}

\subsection{Properties of the potential in our case}

The preceding sections give the following properties for all sufficiently large $t$.

\begin{enumerate}
    \item\label{it:1} $|P(t)-h(t)|\leq C$, and $h(t)\geq st-C$.
    \begin{proof}
        For the first assertion, \eqref{eq:potential comp} gives $-\delta\leq P(t)-h(t)\leq (d-1)c$, so a fixed constant bounds the absolute value of their difference. As for the second assertion, \cref{lem:lower bound for slope} with $r=d$, $\mu=m$ gives $\sigma=s$, and
        \[
        h(t)\geq h(u_0)+s(t-u_0)-C_0.
        \]
        The quantity $h(u_0)-su_0-C_0$ is constant. Thus both assertions hold with
        \[
        C=\max\{(d-1)c,\delta,C_0+su_0-h(u_0)\}.
        \]
    \end{proof}
    \item\label{it:2} If $k$ is active, there exist $k$ independent integer vectors with $\|A_tz\|_\square\leq C'e^{h(t)}$. This is \cref{lem:short ind sys}, where we can choose the constant
    \[
    C'=\max_{1\leq k\leq d}C_k.
    \]
    The resulting system spans the unique minimizing subspace if $k<d$, and spans $\mathbb{R}^d$ in the case $k=d$.
    \item\label{it:3} If rank $k<d$ is active at time $t$, then by \cref{lem:active rank and gap} $f_{k+1}(t)-f_k(t)>2\delta$, and by \cref{lem:min height of rank k}, the subspace attaining the minimum is equal to $U_k(t)$. By continuity the inequality persists throughout some interval $J$, and so the subspace $U_k(t)$ stays constant.
    \item\label{it:4} On an affine piece of $P$ its slope has the form $-n+dl/k$. By \cref{lem:l=0}, $l\geq 1$. If all the active ranks are at most $m$, then this slope is at least $u$.
\end{enumerate}

\subsection{Exclusion of larger ranks}

We prove \cref{problem 2 sol} by contradiction. Assume that $\frac{m^2}{m+1}<\hat{\omega}(\xi)$. Choose a number $\alpha$:
\begin{equation}\label{eq:alpha}
    \frac{m^2}{m+1}<\alpha<\hat{\omega}(\xi).
\end{equation}
Put
\begin{equation}\label{beta}
    \beta=\frac{m-n\alpha}{1+\alpha}<\beta_*.
\end{equation}
The main goal of this subsection is to prove that under this assumption all ranks active at sufficiently large times do not exceed $m+1$.

We shall also use Hadamard's inequality in the following form
\begin{equation}\label{eq:hadamard}
    \|a_1\wedge\dots\wedge a_p\|_2\leq \prod_{j=1}^p \|a_j\|_2.
\end{equation}

\begin{lemma}\label{lem:m+2 rank est}
    There exists a constant $c_\xi>0$ such that for any independent $w_1,\dots,w_{m+2}\in\mathbb{Z}^d$ and all $t\geq0$,
    \begin{equation}\label{eq:m+2 rank est}
        \|A_tw_1\wedge\dots\wedge A_t w_{m+2}\|_2\geq c_\xi e^{\kappa t}.
    \end{equation}
\end{lemma}
\begin{proof}
    Put $p=m+2\leq d$, and $\eta=w_1\wedge\dots\wedge w_p$. Since the $w_i$ have integer coordinates, the coordinates of $\eta$ in the standard exterior basis are integers. Moreover, $\eta\neq0$ because the $w_i$ are linearly independent. Hence
    \[
    \|\eta\|_2\geq 1.
    \]

    Let $B=\bigwedge^p \Phi$. Define the positive number
    \[
    c_\xi=\|B^{-1}\|_2^{-1}>0.
    \]
    Notice that this constant depends only on $\xi$. Using $\eta=B^{-1}(B\eta)$ and the definition of the operator norm gives
    \[
    \|B\eta\|_2\geq \frac{\|\eta\|_2}{\|B^{-1}\|_2}\geq c_\xi.
    \]
    Now decompose $\theta=B\eta$ orthogonally as
    \[
    \theta=\sum_{l=0}^p \theta_l, \quad \theta_l\in\bigwedge^{p-l} E\wedge \bigwedge^l F.
    \]
    Since $\dim E=m$, and $p=m+2$, the degree restriction gives $p-l\leq m$, or $l\geq 2$. On every such component the operator $\bigwedge^p D_t$ multiplies by the scalar
    \[
    e^{(-np+dl)t}.
    \]
    For $t\geq 0,d>0,l\geq2$, this factor is at least $e^{(-np+2d)t}$. The exponent simplifies to
    \[
    -np+2d=-n(m+2)+2(m+n)=m(2-n)=\kappa.
    \]
    Using orthogonality gives
    \[
    \|A_tw_1\wedge\dots\wedge A_tw_p\|_2^2=\sum_l e^{2(-np+dl)t}\|\theta_l\|_2^2\geq e^{2\kappa t}\sum_l \|\theta_l\|_2^2=e^{2\kappa t}\|\theta\|_2^2\geq c_\xi ^2e^{2\kappa t}.
    \]
    Taking the square root proves \eqref{eq:m+2 rank est}.
\end{proof}

\begin{proposition}\label{pr:h est}
    For all sufficiently large $t$
    \begin{equation}\label{eq:h est}
        h(t)\leq\beta t.
    \end{equation}
\end{proposition}
\begin{proof}
    From \eqref{eq:exp formulas} it follows that for all sufficiently large $i$
    \[
    \frac{-\log L_i}{\log X_{i+1}}>\alpha, \quad \log L_i+\alpha\log X_{i+1}<0.
    \]
    Also we have $m-\beta =d\alpha/(1+\alpha)$. Using $h(q_i)=\log L_i+mq_i$ and $dq_i=\log X_{i+1}-\log L_i$, we obtain
    \[
    h(q_i)-\beta q_i=\log L_i+\frac{d\alpha}{1+\alpha}q_i=\frac{\log L_i+\alpha\log X_{i+1}}{1+\alpha}<0.
    \]
    Thus the estimate holds at every sufficiently late point $q_i$.

    Now it suffices to prove it on $[q_i,q_{i+1}]$. \eqref{h1 formula} tells us that on this interval
    \[
    h(v)=\max(\log X_{i+1}-nv,\log L_{i+1}+mv).
    \]
    It is the maximum of two affine functions. Therefore, it is a convex, piecewise affine function. The same is true for $g(v)=h(v)-\beta v$. Every point on the considered interval has the form $v=\theta q_i+(1-\theta)q_{i+1},0\leq\theta\leq1$. Convexity gives
    \[
    g(v)\leq \theta g(q_i)+(1-\theta)g(q_{i+1})\leq 0,
    \]
    since $g(q_i),g(q_{i+1})<0$, as shown before.
\end{proof}

\begin{lemma}
    Under \eqref{eq:alpha} all sufficiently late active ranks are at most $m+1$.
\end{lemma}
\begin{proof}
    Suppose that rank $k\geq m+2$ is active at a sufficiently late time $t\geq0$. By \cref{it:2}, we can obtain $k$ independent integer vectors $w_1,\dots,w_k$, for which $\|A_tw_j\|_\square\leq C'e^{h(t)}$. Choose $p=m+2$ members of this system. For them \eqref{block and eucl norm comp} gives
    \[
    \|A_tw_j\|_2\leq \sqrt{2}\|A_tw_j\|_\square\leq\sqrt{2}C'e^{h(t)}.
    \]
    Hadamard's inequality \eqref{eq:hadamard} therefore gives
    \[
    \|A_tw_1\wedge\dots\wedge A_tw_p\|_2\leq (\sqrt{2}C')^pe^{p h(t)}.
    \]

    This and \cref{lem:m+2 rank est} give
    \[
    c_\xi e^{\kappa t}\leq (\sqrt{2}C')^pe^{p h(t)}.
    \]
    Taking the logarithm gives
    \[
    h(t)\geq\frac{\kappa}{p}t-D,\quad D=\max\left( 0,\log(\sqrt{2}C')-\frac{\log c_\xi}{p} \right).
    \]

    Recall the estimate $h(t)\leq\beta t$ from \cref{pr:h est}. Combining it with the last inequality gives us
    \[
    \left( \frac{\kappa}{m+2}-\beta \right)t\leq D.
    \]
    However \eqref{eq:const relations} and \eqref{beta} show that the coefficient on the left is positive. Hence this inequality is impossible for sufficiently large $t$.
\end{proof}

\subsection{Transitions between subspaces}

Choose a number $t_0\geq0$ such that for all $t\geq t_0$ all aforementioned properties and lemmas hold. In particular, no active rank exceeds $m+1$. Define the set
\begin{equation}\label{eq:A set}
    \mathcal{A}=\{t\geq t_0: \text{rank}\ m+1 \ \text{is active at time}\ t\}.
\end{equation}
For each $t\in \mathcal{A}$ we denote by $Q(t)$ the corresponding rational subspace. Its uniqueness follows from \cref{it:3}, since $m+1<d$.

\begin{lemma}\label{lem:transition}
    There exists a sequence of intervals $[a_j,b_j]$ such that $a_j\to\infty$ and the following properties hold:
    \begin{equation}\label{eq:int properties}
        a_j,b_j\in\mathcal{A},\quad a_j<b_j,\quad (a_j,b_j)\cap\mathcal{A}=\emptyset,\quad Q(a_j)\neq Q(b_j).
    \end{equation}
    Moreover, there exists a constant $C>0$ such that on each such interval
    \begin{equation}\label{eq:int ineq}
        h(b_j)-h(a_j)\geq u(b_j-a_j)-C.
    \end{equation}
\end{lemma}
\begin{proof}
    Consider the function
    \[
    F(t)=\frac{G_{m+1}(t)}{m+1}+(d-m-1)c-P(t).
    \]
    It is always nonnegative by definition of the potential. It equals zero precisely at times $t$ when rank $m+1$ is active. Moreover, $F$ is continuous and locally piecewise affine, since the functions $G_{m+1}$ and $P$ are continuous and piecewise affine. Consequently, $\mathcal{A}=F^{-1}\{0\}\cap [t_0,+\infty)$ is closed.

    Consider an interval on which $F$ is affine. Its set of zeros on this interval is either empty, a single point, or the whole interval. Therefore, $\mathcal{A}$ is a union of closed intervals and isolated points. Consider a connected component of $\mathcal{A}$. By \cref{it:3}, the space $Q(t)$ is constant in some neighborhood of $t$ in this connected component. A locally constant function on a connected space is constant. Hence $Q(t)$ is constant on every connected piece of $\mathcal{A}$.

    Now we show that $\mathcal{A}$ is unbounded. Suppose that it is bounded. Choose $T_1\geq t_0$ such that no active rank at time $t\geq T_1$ exceeds $m$. By \cref{it:4}, the slope of $P$ on every affine piece is at least $u$. For $t\geq T_1$ this gives
    \[
    P(t)\geq P(T_1)+u(t-T_1).
    \]
    And \cref{it:1} gives
    \[
    h(t)\geq ut-D_1,
    \]
    where $D_1$ is a fixed constant. This along with the estimate $h(t)\leq\beta t$ implies
    \[
    (u-\beta)t\leq D_1.
    \]
    However $u>\beta$. Thus, this is impossible.

    Now we show that the subspace $Q(t)$ cannot be eventually constant for $t\in\mathcal{A}$. Suppose, to the contrary, that $Q(t)=Q$ for all $t\geq T_2$ in $\mathcal{A}$. Then $\dim Q=m+1<d$, hence it is a proper subspace, and $\mathcal{L}_\xi\not\subset Q$. Therefore, $\rho=\dim(\Phi Q\cap E)=\dim(Q\cap\mathcal{L}_\xi)\leq m-1$. By \cref{lem:largest err deg}, the largest index $l$ for which $K_l(Q)\neq 0$ equals $l_*=m+1-\rho\geq 2$.
    
    Similarly to the proof of \cref{tail theorem}, for this $Q$, the ratios of smaller degree terms to $l_*-\text{degree}$ term tend to $0$. Thus after a time $t_Q$ for all $t\in\mathcal{A}$, the slope of $P$ on an affine piece is
    \[
    -n+\frac{dl_*}{m+1}\geq-n+\frac{2d}{m+1}\geq -n+\frac{d}{m}=u.
    \]
    If a rank $k\leq m$ is active, then the slope of $P$ is at least $u$ by \cref{it:4}. Now we arrive at a contradiction, similarly to the previous step.

    Consequently, there are arbitrarily large $t$ that do not belong to $\mathcal{A}$. Otherwise, some half-line would belong to $\mathcal{A}$. It is contained in a single connected component of $\mathcal A$, so $Q(t)$ is constant on that half-line.

    Since $F$ is locally piecewise affine, only finitely many connected components of $\mathcal A$ meet any compact interval. Now list all the connected components of $\mathcal{A}$ in their order on the number line $C_1,C_2,\dots$. There are infinitely many of them, since $\mathcal{A}$ is unbounded and contains no half-line.

    Choose $t\in C_i$, and put $Q_i=Q(t)$. As shown before, this subspace is the same regardless of the choice of $t$. Since $Q(t)$ is not eventually constant on $\mathcal{A}$, there are infinitely many indices $j$ such that $Q_j\neq Q_{j+1}$. For every such index put
    \[
    a_j=\max C_j,\quad b_j=\min C_{j+1}.
    \]
    These endpoints belong to $\mathcal{A}$, since it is closed. They satisfy $a_j<b_j$, since the components are different, and $(a_j,b_j)\cap\mathcal{A}=\emptyset$, since the components are in order. This proves all the properties in \eqref{eq:int properties}.

    Now we prove \eqref{eq:int ineq}. Since the active rank on $(a_j,b_j)$ does not exceed $m$, by \cref{it:4} the slope of $P$ on every affine piece in this interval is at least $u$. Hence,
    \[
    P(b_j)-P(a_j)\geq u(b_j-a_j).
    \]
    \eqref{eq:potential comp} gives
    \[
    h(b_j)-h(a_j)\geq P(b_j)-P(a_j)-(d-1)c-\delta\geq u(b_j-a_j)-C,\quad C=(d-1)c+\delta,
    \]
    as required.
\end{proof}

\begin{lemma}
    For every interval $[a_j,b_j]$ supplied by \cref{lem:transition}, we have
    \begin{equation}\label{eq:transition ineq}
        (m+1)h(a_j)+h(b_j)+n(b_j-a_j)\geq\kappa a_j-C,
    \end{equation}
    where the constant does not depend on $j$.
\end{lemma}
\begin{proof}
    For convenience, put $a=a_j,b=b_j,\tau=b_j-a_j>0$. \Cref{it:2} supplies $m+1$ independent integer vectors $w_1,\dots,w_{m+1}$, spanning $Q(a)$, that satisfy
    \[
    \|A_aw_i\|_\square\leq C'e^{h(a)}.
    \]
    Applying \cref{it:2} at time $b$, choose $m+1$ independent integer vectors spanning $Q(b)$ and satisfying $\|A_bz\|_\square\leq C'e^{h(b)}$. Since $Q(a)\neq Q(b)$, one of these vectors, say $v$, lies outside $Q(a)$. Hence $w_1,\dots,w_{m+1},v$ are linearly independent.

    Consider
    \[
    A_aA_b^{-1}=D_a\Phi\Phi^{-1}D_{-b}=D_{a-b}=D_{-\tau}.
    \]
    If $A_bv=(u_b,v_b)$, then $A_av=(e^{n\tau}u_b,e^{-m\tau}v_b)$. Since $\tau>0,m,n>0$, both scale factors are at most $e^{n\tau}$. Hence,
    \[
    \|A_av\|_\square=\max(e^{n\tau}\|u_b\|_2,e^{-m\tau}\|v_b\|_2)\leq e^{n\tau}\|A_bv\|_\square\leq C'e^{h(b)+n\tau}.
    \]
    All $m+2$ vectors can now be estimated at the same time $a$.

    Applying \cref{lem:m+2 rank est}, Hadamard's inequality, and \eqref{block and eucl norm comp} gives
    \[
    c_\xi e^{\kappa a}\leq \|A_aw_1\wedge\dots\wedge A_aw_{m+1}\wedge A_a v\|_2\leq (\sqrt{2}C')^{m+2}e^{(m+1)h(a)+h(b)+n\tau}.
    \]
    Taking the logarithm and rearranging gives
    \[
    (m+1)h(a)+h(b)+n\tau\geq\kappa a-C,\quad C=\max(0,(m+2)\log(\sqrt{2}C')-\log c_\xi),
    \]
    which is exactly \eqref{eq:transition ineq}.
\end{proof}

\subsection{End of the proof}

\begin{proof}[Proof of \Cref{problem 2 sol}]
    Choose an interval $[a,b]$ from \cref{lem:transition}. Set $x=h(a),y=h(b),\tau=b-a$. We have three uniform estimates
    \begin{equation}\label{eq:est 1}
        x-sa\geq -C,
    \end{equation}
    \begin{equation}\label{eq:est 2}
        y-x-u\tau\geq -C,
    \end{equation}
    \begin{equation}\label{eq:est 3}
        (m+1)x+y+n\tau-\kappa a\geq -C.
    \end{equation}
    The first comes from \cref{it:1}, the second is \eqref{eq:int ineq}, and the third is \eqref{eq:transition ineq}.

    Multiply these inequalities respectively by
    \[
    m^2-1,\quad m(m+1),\quad 1.
    \]
    Adding the resulting inequalities eliminates $x$ and yields
    \[
    \Delta_m(y-\beta_*(a+\tau))\geq -((m^2-1)+m(m+1)+1)C,
    \]
    or
    \[
    h(b)\geq \beta_*b-C.
    \]
    Combine this with \eqref{eq:h est} to obtain
    \[
    (\beta_*-\beta)b_j\leq C.
    \]
    However, $\beta_*-\beta>0$, and $b_j$ can be chosen arbitrarily large. Hence this inequality is impossible for $b_j$ large enough. This contradicts our assumption \eqref{eq:alpha}. Thus
    \[
    \hat{\omega}(\xi)\leq\frac{m^2}{m+1}.
    \]
\end{proof}

\begin{corollary}\label{cor:C=BC uniform exp}
    For $m,n\geq 2$
    \begin{equation}\label{eq:C+BC uniform exp}
        \Sigma_{\hat{\omega}}(\mathcal{C}_{m,n})=\Sigma_{\hat{\omega}}(\mathcal{BC}_{m,n})=\left[ \frac{m}{n},\frac{m^2}{m+1} \right].
    \end{equation}
\end{corollary}
\begin{proof}
    By \cref{problem 2 sol}, if condition \eqref{cond C} holds for totally irrational $\xi$, then its uniform exponent lies in this interval. The same is obviously true if both conditions \eqref{cond B} and \eqref{cond C} hold. In \cite[Theorem~2.9(ii)]{MMS} it is proven that any uniform exponent from this interval is admissible when both conditions hold. This proves \eqref{eq:C+BC uniform exp}.
\end{proof}

\section{\texorpdfstring
  {Proof of \Cref{problem 3 sol}}
  {Proof of Theorem \ref*{problem 3 sol}}}

\begin{proof}[Proof of \Cref{problem 3 sol}]
    The first row follows from \cref{C=BC ord exp,cor:C=BC uniform exp}.

    For the second row, \cite[Theorem~2.2]{MMS} states that \eqref{cond C} implies \eqref{cond B}, and in this case $\xi$ is badly approximable. In particular, $\mathcal{C}_{m,1}=\mathcal{BC}_{m,1}$, and $\omega(\xi)=\hat{\omega}(\xi)=m$.

    As for the final row, \cite[Theorem~2.2]{MMS} gives bad approximability if \eqref{cond B} holds. Hence in this case $\omega(\xi)=\hat{\omega}(\xi)=\frac{1}{n}$. \Cref{C=BC ord exp} gives
    \[
    \Sigma_\omega (\mathcal{C}_{1,n})=\left[ \frac{1}{n}, 1\right].
    \]
    Finally, for the uniform exponent in the case $m=1,n\geq2$, with condition \eqref{cond C}, Moshchevitin and Akhunzhanov proved in \cite[Proposition~2]{AkhunzhanovMoshchevitin} that $\hat{\omega}(\xi)\leq 1/2$. Later this case was fully solved by Summerer. In \cite{Summerer} it is proven that 
    \[
    \Sigma_{\hat{\omega}}(\mathcal{C}_{1,n})=\left[ \frac{1}{n},\frac{1}{2} \right].
    \]
    This completes the proof.
\end{proof}

Therefore, when $m,n\geq2$, the extra condition \eqref{cond B} does not change the spectra. In the case $n=1$ the situation is the same. However, when $m=1,n\geq2$, it changes the ordinary spectrum. The uniform spectrum changes only when $n\geq3$. These statements concern the individual spectra, not the joint spectrum of $(\omega,\widehat{\omega})$.

\section*{AI assistance}

The list of problems in the original paper was given by the author to OpenAI's ChatGPT-5.6 with the goal of selecting which problem to work on. Its first choice was Problem~2. After a few prompts the model was able to devise a solution to this problem in two special cases. A key ingredient was a weaker version of \cref{th:block lemma}.

The same problems were then provided to ChatGPT-6 Astra, along with the aforementioned solutions. The model was able to develop solutions to all three problems, building upon and generalizing the work done before.

All the proofs, including their logic and calculations, were carefully checked by the author. The initial text was rewritten to improve clarity and readability. The author takes full responsibility for the correctness of the results of this paper.

\end{document}